\documentclass[11pt]{amsart}
\usepackage{xypic}
\usepackage{amssymb}   
\usepackage{amsmath}
\usepackage{amsthm}
\usepackage{times}

\theoremstyle{theorem}
\newtheorem{Theorem}{Theorem}[section]

\newtheorem{Proposition}[Theorem]{Proposition}
\newtheorem{Corollary}[Theorem]{Corollary}

 \theoremstyle{definition}
\newtheorem{Definition}[Theorem]{Definition}

\numberwithin{equation}{section}

\newcommand{\Pic}{\operatorname{Pic}}
\newcommand{\Bl}{\operatorname{Bl}}

\renewcommand{\O}{{\mathcal O}}

\newcommand{\Z}{{\mathbb Z}}

\newcommand{\C}{{\mathbb C}}
\newcommand{\p}{{\mathbb P}}

\newcommand{\codim}{\operatorname{codim}}

\newcommand{\map}{\dasharrow}

\newcommand{\Chow}{\operatorname{Chow}}
\newcommand{\Num}{\operatorname{Num}}

\newcommand{\Supp}{\operatorname{Supp}}

\def\leq{\leqslant}
\def\geq{\geqslant}
\def\deft{~}

\def\bibaut#1{{\sc #1}}

\def\phi{\varphi}

\begin{document}

\title[CC-manifolds]{Conic-connected manifolds}
\author[Paltin Ionescu]{Paltin Ionescu*}
\address{University of Bucharest, Faculty of Mathematics and Computer Science, 14 Academiei str., RO-010014 Bucharest, and  Institute of Mathematics of the Romanian Academy, P.O. Box 1-764, RO-014700 Bucharest, Romania}
\email{Paltin.Ionescu@imar.ro}
\author[Francesco Russo]{Francesco Russo**}
\thanks{*Partially  supported  by the Italian Programme ``Incentivazione alla mobilit\`{a} di studiosi stranieri e italiani residenti all'estero"}
\address{Dipartimento di Matematica e Informatica\\
Universit\` a degli Studi di Catania\\
Viale A. Doria, 6\\
95125 Catania\\ Italy
}

\thanks{**Partially  supported  by CNPq
(Centro Nacional de Pesquisa), Grants 300961/2003-0, 308745/2006-0 and 474475/2006-9, and by
PRONEX/FAPERJ--Algebra Comutativa e Geometria Algebrica}
\email{frusso@dmi.unict.it}
\begin{abstract}
We study a particular class of rationally connected manifolds, $X\subset \p^N$, such that two general points $x,x' \in X$ may be joined by a conic contained in $X$. We prove that these manifolds are Fano, with $b_2\leq 2$. Moreover, a precise classification is obtained for $b_2=2$. Complete intersections of high dimension with respect to their multi-degree provide examples for the case $b_2=1$. The proof of the classification result uses a general characterization of rationality, in terms of suitable covering families of rational curves.   
\end{abstract}

\maketitle
\thispagestyle{empty}
\begin{center}\begin{minipage}[l]{10.5cm}
\hskip12pt {\it \footnotesize A Gaetano Scorza che un secolo fa 
aveva colto l'importanza delle variet\`{a} razionalmente connesse
considerando  la classe particolare di variet\`{a} conicamente
connesse formata da quelle di {\it ultima specie}.

\hskip12pt "Invece per le $V_4$ di prima e terza specie arrivo a caratterizzarle tutte
valendomi della teoria dei sistemi lineari sopra una variet\`{a} algebrica e, per le ultime,
della circostanza che esse contengono un sistema $\infty^6$ di coniche cos\`{i} che 
per ogni loro coppia di punti passa una e una sola conica.
Inoltre la natura dei ragionamenti \`{e} tale da mostrare come i
risultati ottenuti possano estendersi, almeno per la maggior parte,
alle variet\`{a} (di {\it prima} e {\it ultima} specie) a un numero qualunque di
dimensioni, ...", \cite[{\it Opere Scelte},\deft vol.\deft 1,\deft p.\deft253]{Scorza2}.}
\end{minipage}\end{center}
\normalsize

\vskip10pt

\section*{Introduction}

\thispagestyle{empty}

A projective manifold $X$ is rationally connected if any two given (general) points of it may be joined by a rational curve. Since their introduction by Koll\'ar, Miyaoka, and Mori \cite{KMM}, rationally connected manifolds were intensively studied, as they naturally generalize both classical notions of unirational and Fano manifold. When such an $X$ is embedded in some projective space, a measure of its complexity is given by the minimal degree of a family of rational curves possessing the above defining property. As a first instance we have linear spaces, characterized by the presence of a line through any two distinct points. The next case is what we call {\it conic-connected} manifolds: They are embedded manifolds $X\subset \p^N$ such that for any two general points $x, x' \in  X$, there is some irreducible conic joining $x,x'$ and contained in $X$.

The main results of this paper, Theorem\deft\ref{delta1} and Corollary~\ref{C4.2}, show that, over $\C$, conic-connected manifolds  are Fano, with second Betti number $b_2\leq 2$. 
Moreover, those with $b_2=2$ are completely described in the following short list:
Segre products of two projective spaces and  their hyperplane sections, or the~inner projections, from a linear space, of the Veronese variety $v_2(\p^n)$. When $b_2=1$, unless $X \simeq v_2(\p^n)$, the Picard group is generated by the hyperplane section and the index of $X$ is at least $\frac{\dim(X)+1}2$. 
Conversely, from a recent result by Bonavero--H\" oring, see \cite{BH}, it follows that any $n$-dimensional complete intersection which is Fano, of index at least $\frac{n+1}2$, is conic-connected; see Proposition 2.4. Moreover, the work of Hwang--Kebekus, see \cite{HK}, shows that any Fano $n$-fold of index greater than $\frac{2n}3$ is conic-connected; see Proposition 2.5. Applications of Theorem\deft\ref{delta1} to the classification of secant and dual defective manifolds, special Cremona transformations, etc. may be found in \cite{Ru, IR}. 
 
The proof of Theorem\deft\ref{delta1} naturally fits into the context of Mori Theory, as established by the ground-breaking papers \cite{Mori, Mori1}. 
Besides the use of such powerful tools as the Cone Theorem, contractions of extremal rays, etc., we also need the very efficient method of algebraizing the equivalence relation associated to a covering family of curves. 
This technique, pioneered by Campana \cite{Ca}, was further developed and successfully applied by Koll\' ar, Miyaoka, and Mori in their seminal paper \cite{KMM}. 
The recent characterization of projective space, belonging to the same circle of ideas, and due to Cho, Miyaoka, and Shepherd-Barron (see \cite{CMSB}), is also used in the proof of Theorem\deft\ref{delta1}. 
Finally, a key ingredient of the argument is our Theorem\deft\ref{4.5}, a ``global version"  of Theorem\deft\ref{4.4}. 

Theorem\deft\ref{4.4} which, we hope, has some independent interest, generalizes \cite[Theorem\deft4.2]{IN} (see also \cite{Io2}).  
Its content is a criterion of rationality: 
$X$ is rational  if and only if it admits a covering family of rational $1$-cycles, all passing through $x$, all smooth at $x$ and such that the general cycle of the family is a curve, uniquely determined by its tangent line at $x$. 
In Theorem\deft\ref{4.5}, under suitable assumptions, the birational isomorphism  between $X$ and $\p^n$ is shown to be  biregular along a general (smooth) curve of the given family; moreover, 
the general member of the family corresponds to a line in $\p^n$.

 We acknowledge our intellectual debt to the pioneering work of G.\ Scorza \cite{Scorza1, Scorza2}, as well as to all the previous
papers mentioned above.    
        
We work over $\C$; note however that Theorem\deft\ref{4.4} holds in any characteristic. 

\section{A rationality criterion}\label{covering}

Let $X$ be a projective variety. A family of rational $1$-cycles on $X$, $\mathcal{F} \to \mathcal{C}$, is given by a closed subset $\mathcal{C}\subset {\rm Chow}(X)$ such that for each $1$-cycle of the family, all its irreducible components are rational. Moreover, we always assume that $\mathcal{C}$ is {\it irreducible} and the {\it general member} of $\mathcal{C}$ is an {\it irreducible reduced} (rational) curve. 
  
\begin{Definition}
A {\it covering family} (of rational $1$-cycles, as above) on $X$ is a family  $\mathcal{F} \to \mathcal{C}$ as above, such that the tautological morphism $\mathcal{F} \to X$ is surjective.  
\end{Definition}

Let $x\in X$ be a {\it smooth} point.
\begin{Definition}
An {\it $x$-covering family} on $X$ is a covering family $\mathcal{F} \to \mathcal{C}$ as above, such that, moreover, $x\in \Supp (C)$ for any $C\in \mathcal{C}$ and the general member of $\mathcal{C}$ is smooth at $x$. 

An $x$-covering family is called {\it maximal} if $\mathcal{C}$ is an irreducible component of the closed subscheme of ${\rm Chow}_x(X)$, whose points correspond to rational $1$-cycles.

An $x$-covering family is said to be {\it smooth} if all its $1$-cycles are smooth {\it at $x$}. 

An $x$-covering family is said to satisfy the {\it infinitesimal uniqueness property} if the general member of the family is uniquely determined by its tangent space at $x$.
\end{Definition}

Recall that $X$ is {\it rationally connected} if through any two (general) points  of $X$ there passes a rational curve. 
Saying that $X$ is rationally connected is equivalent to the existence of an $x$-covering family.
One can also see that the existence of an $x$-covering family satisfying the infinitesimal uniqueness property is equivalent to $X$ being {\it unirational}.
The next result gives a criterion of {\it rationality}, in the spirit of the papers \cite{IN, Io2}. It may be seen as a birational version of a series of characterizations of $\p^n$ in terms of families of rational curves, started in Mori's famous paper \cite{Mori}.     

\begin{Theorem}\label{4.4}
Let $X$ be a projective variety. The following conditions are equivalent: 
\begin{itemize}
\item[(i)] $X$ is rational;

\item[(ii)] for some $x \in X$, $X$ admits a smooth $x$-covering family satisfying the infinitesimal uniqueness property.
\end{itemize}  
\end{Theorem}

\begin{proof}
To see that (i) implies (ii), note that a neighbourhood of a general point $x\in X$ is isomorphic to some neighbourhood of a point $y\in \p^n$. Simply consider the family induced on $X$ by the lines of $\p^n$ passing through $y$. 

(ii) implies (i), the non-trivial implication, is proved in the following three steps.
\vskip4pt

\noindent{\it Step} I. The basic diagram. 

 Let $\pi: \mathcal{F} \to \mathcal{C}$ be our $x$-covering family and let $\phi: \mathcal{F} \to X$ be the tautological morphism.

Note that $\pi$ admits a section ${\mathcal E}$, which is contracted by $\phi$ to the point $x$.  Consider the blowing-up $\sigma: X'\to X$ of $X$ at $x$. Via \cite[Lemma\deft4.3]{IN}, the universal property of the blowing-up shows that we have a morphism $\psi : {\mathcal F} \to X'$, such that $\sigma \circ \psi = \phi$. So we have the following diagram
\begin{equation*}\label{joindiagram1}\raisebox{.7cm}{\xymatrix{
&\mathcal{F}  \ar[d]_\pi \ar[dr]^\phi\ar[r]^\psi&X'\ar[d]^\sigma\\
&\mathcal{C}&X.
}}
\end{equation*}
In particular, $\psi $ maps the section ${\mathcal E}$   to $E$, the exceptional locus of $\sigma$. Let $\psi _0 : {\mathcal E} \to E$ be the restriction of $\psi $. 

\vskip4pt

\noindent{\it Step} II. $\psi _0$ is generically finite and surjective; in particular, $\dim ({\mathcal C})= n-1$.

This follows from an elementary, but very useful remark, due to Kebekus \cite{Kebekus}, proof of Theorem\deft3.4. Indeed, consider a general point $p\in E\simeq \p ( T ^*_x (X))$ and let $C\subset {\mathcal E} \simeq {\mathcal C}$ be a curve such that $\psi _0(C)= p$. Take $\widetilde C$, the normalization of $C$, and consider the surface $S$ over $\widetilde C$ which is got by base-change from our family. Let $D\subset S$ be the induced section for $\pi_S : S\to \widetilde C$ and note that $S$ is smooth along $D$ by our hypothesis. The restriction of the tangent morphism of $\phi_S: S\to X$ induces a morphism $T_{\phi_S}: N_{D|S} \to l_p\simeq \C$, where $N_{D|S}$ is the (geometric) normal bundle of $D$ in $S$ and $l_p$ is the line in $ T _x (X)$ corresponding to $p\in E$. Since $N_{D|S}$ is not trivial, the above map has a zero. The corresponding curve of the family is singular at $x$. This is a contradiction. To conclude, note also that $\dim ({\mathcal F}) \geq n$, so $\dim ({\mathcal E})\geq n-1$, while $ \dim (E) = n-1$. In particular, $\dim ({\mathcal C})= n-1$ and $\phi$ (or $\psi $) is generically finite.

\vskip4pt
\noindent{\it Step} III. 
Conclusion.

Note first that the infinitesimal uniqueness property of our family translates into saying that the map $\psi_0$ is birational. 
Observe also that we have ${\psi }^{-1} (E) = {\mathcal E} $. Now, use properness of $\psi $ and the fact that $\psi _0$ is birational to deduce that $\psi $ is birational. We have already seen that ${\mathcal C}$ is birational to  $E \simeq \p^{n-1}$.
Replacing $\mathcal {F}$ by its normalization $\widetilde {\mathcal F}$ and using the section ${\mathcal E}$, we find that $\widetilde {\mathcal F}$ is generically a $\p^1$-bundle over $\mathcal{C}$. Thus $X$ is rational, being birational to $\widetilde {\mathcal F}$. 
\end{proof}

The next result is a ``more global" version of Theorem\deft\ref{4.4}; it improves \cite[Theorem\deft4.2]{IN}. In the terminology of \cite{IV}, (ii) says that the model $(X,C)$ is equivalent to $(\p^n, {\rm line})$. Let us recall the following definition from \cite{BaBeIo}:

\begin{Definition}[\cite{BaBeIo}]\label{bbi} A smooth rational curve $C \subset X$, where $X$ is a projective manifold of dimension $n$, is a {\it quasi-line} if $N_{C|X}\simeq \bigoplus^{n-1}_1 \mathcal{O}_{\p^1}(1).$
\end{Definition}

\begin{Theorem}\label{4.5}
Let $X$ be a projective manifold and let $\mathcal{C}$ be a smooth maximal $x$-covering family of it. Let $\mathcal{D} \subset \mathcal{C}$ denote the family of reducible cycles in $\mathcal {C}$. 
\begin{itemize}
\item[(i)] If the general member of $\mathcal{C}$ is smooth, it is a quasi-line (this is always the case if $\dim (X) \geq 3$).
\item[(ii)] Assume moreover that the general member of $\mathcal{C}$ meets the locus of $\mathcal {D}$ only at $x$ (this is true, e.g.\ if $\codim _{\mathcal{C}} (\mathcal{D}) \geq 2$). Then, for any quasi-line $C\in \mathcal{C}$, 
 there are a neighbourhood $U$ of $C$ in $X$ and a neighbourhood $V$ of a line $l\subset \p^n$, together with an isomorphism $U\simeq V$, sending~$C$\deft  to\deft $l$.  
\end{itemize}
\end{Theorem}
\begin{proof} We use the setting and notation from the proof of Theorem\deft\ref{4.4}.

(i) Let $C\in \mathcal{C}$ be a general smooth curve. As $\mathcal{C}$ is a maximal $x$-covering family, $N_{C|X}$ is ample. In Step II from the proof of Theorem\deft\ref{4.4} we saw that $\dim (\mathcal{C})=n-1$, so $C$ is a quasi-line.

(ii) We first prove the following:

  {\it Claim}: There is a unique quasi-line from our given family passing through a general point of $X$ (i.e. the map $\phi$ is birational).

  Let $D\subset X'$ be the locus of points belonging to the images via $\psi$ of reducible cycles from $\mathcal{C}$. By (i), a general curve $C$ from our family $\mathcal{C}$ is a quasi-line in $X$. By the assumption in (ii),  if $\tilde C$ denotes its proper transform on $X'$,   we have $\tilde {C} \cap D = \emptyset$ for $C$ general.  So any point $c\in \pi (\psi^{-1}(\tilde {C}))$ represents an irreducible member of the family $\mathcal{C}$.
Now, the birationality of the map $\phi$ follows from the proof of \cite[V.3.7.5.]{Kollar}, a result due to Miyaoka. So the Claim is proved.

  The Claim allows us to apply \cite[Theorem\deft2.5]{Io2} in order to get the conclusion in (ii).
Finally, observe that $\codim _{\mathcal{C}} (\mathcal{D}) \geq 2$ implies $\codim _{X'} (D) \geq 2$. Moreover, if $C$ is a general member of $\mathcal{C}$, the normal bundle of $\tilde{C}$ in $X'$ is trivial; see, for instance, \cite[Lemma 1.6]{IN}. So, by  \cite[II.3.7]{Kollar}, we have $\tilde {C} \cap D = \emptyset$ for $C$ general. 
\end{proof}

\section{Classification of conic-connected manifolds}

\begin{Definition}[cf.\ also \cite{KS}] A smooth  irreducible non-degenerate projective
variety $X\subset\p^N$ is said to be a {\it conic-connected
manifold}, briefly a {\it CC-manifold}, if through two general
points of $X$ there passes an irreducible conic contained in $X$. 
\end{Definition}

The following classification theorem is our main result; as CC-manifolds   
are stable under isomorphic projection, we may assume $X$ to be linearly normal.

\begin{Theorem}\label{delta1}
Let $X\subset\p^N$ be a smooth irreducible linearly normal non-degenerate
CC-manifold of dimension $n$. Then either $X\subset\p^N$ is a Fano manifold with
$\Pic(X)\simeq\mathbb{Z}\langle \O_X(1)\rangle $ and of index
$i(X)\geq \frac{n+1}{2}$, or it is
projectively equivalent to one of the following:
\begin{enumerate}
\item[(i)] $\nu_2(\p^n)\subset\p^{\frac{n(n+3)}{2}}$.

\item[(ii)] 
The projection of $\nu_2(\p^n)$ from the linear
space $\langle \nu_2(\p^s)\rangle $, where $\p^s\subset\p^n$ is a linear subspace; equivalently
$X\simeq \Bl_{\p^s}(\p^n)$ embedded in $\p^N$ by the linear system
of quadric hypersurfaces of $\p^n$ passing through $\p^s$;
alternatively $X\simeq\p_{\p^r}(\mathcal{E})$ with
$\mathcal{E}\simeq\O_{\p^r}(1)^{\oplus n-r}\oplus\O_{\p^r}(2)$,
$r=1,2, \ldots, n-1$, embedded by $|\O_{\p(\mathcal{E})}(1)|$. Here $N=\frac{n(n+3)}{2}-\binom{s+2}{2}$ and  $s$ is an integer such that $0\leq s\leq n-2$.

\item[(iii)]  A hyperplane section of the Segre embedding 
$\p^a\times\p^{b}\subset\p^{N+1}$. Here $n\geq 3$ and  $N=ab+a+b-1$, where $a\geq 2$ and $b\geq 2$ are    such that $a+b=n+1$.

\item[(iv)] $\p^a\times\p^b\subset\p^{ab+a+b}$ Segre embedded, where $a,b$ are positive integers such that
$a+b=n$.
\end{enumerate}
\end{Theorem}

\begin{Corollary} \label{C4.2} CC-manifolds are Fano and have second Betti number $b_2\leq 2$;
 when $b_2=2$, they are also rational.
\end{Corollary}

\begin{proof}[Proof of Theorem~{\rm \ref{delta1}}] 
Fix a general point $x\in X$. By definition of CC-manifold, there exists a (maximal) $x$-covering family of conics on $X$; we fix one such family, which we denote by $\mathcal{C}_x$. 

Assume first that there exist lines $l_1$ and $l_2$ 
such that $l_1\cup l_2$ belongs to $\mathcal{C}_x$.

Let $[C]$ denote the numerical
class in $\Num(X)$ of the curve $C$. Thus, if we let [Q] be the class of a general conic from the family $\mathcal{C}_x$, we have $[Q] = [l_1] + [l_2]$.
Let
$V_i$, $i=1,2$, be the irreducible proper component of
$\Chow_{1,1}(X)$ to which $l_i$ belongs, notation as in
\cite{Kollar}.  

As in the first section, let us call the family $V_i$ covering if its curves fill up $X$. Observe that, by the generality of the point $x\in X$, at least one of the families $V_1, V_2$ is so.
 
For each covering family $V_i$, define a relation
$R_i$ on $X$ by saying that two points $x_1,x_2\in X$ are
equivalent if there exists a connected chain of lines from the family $ V_i$ joining the points
$x_1$ and $x_2$. This is the $\langle (U_i)_1\rangle $-relation, $i=1,2,$ for the
normal form associated to the proper connected universal family
$U_i\to V_i$ \cite[Sections IV.4, IV.4.4 and IV.4.8.3]{Kollar}. 

Theorem IV.4.16 of \cite{Kollar} (see also \cite[Theorem
5.9]{Debarre}) yields the existence of an open subset
$X_i^0\subseteq X$ and of a proper morphism with connected fibers
$\pi_i:X^0_i\to Z^0_i$ such~that the $\langle (U_i)_1\rangle $-relation restricts
to an equivalence relation on $X^0_i$ and such that, for every $z\in
Z^0_i$, $\pi_i^{-1}(z)$ coincides with the equivalence relation
class defined by $V_i$, $i=1,2$. Since the family is covering the equivalence
class of a point in $X_i^0$ does not reduce to the point itself.
 Therefore the equivalence class of a general point $x\in X$ is scheme theoretically a
smooth irreducible positive dimensional projective algebraic variety $F^x_i=\pi_i^{-1}(\pi_i(x))$.
Let $a_i=\deg N_{l_i|X}$. Since we are assuming $V_i$ covering, $a_i\geq 0$; the locus of lines (from the family determined by $l_i$) passing through a general point $x$ will be denoted by $C^x_i$. It is a cone of dimension $a_i+1\geq 1$. Clearly, we have $C^x_i \subseteq F^x_i$, so $\dim (F^x_i) \geq a_i+1$ and equality holds if and only if $F^x_i$ is a~linear~space.

\vskip4pt

\noindent{\it Case} I. 
Assume that $[l_1]= [l_2]$. 
 
Let $[l]=[l_i]$, $i=1,2$, and let $Q\in \mathcal{C}_x$ be a general
conic.  Let $R_x$ be the locus of points on $X$ which can be joined
to $x$ by a connected chain of lines whose numerical class is $[l]$.
If $R_x\subsetneq X$, reasoning as above, we could construct an open
subset $X^0\subseteq X$ and a proper morphism with connected fibers
$\pi:X^0\to Z^0$ 
 such that, for every $z\in Z^0$, $\pi^{-1}(z)$ coincides
with the equivalence relation class defined by  the irreducible
proper component of $\Chow_{1,1}(X)$ to which $l$ belongs. From this
it would follow the existence of an effective divisor $D$ passing through a
general point  $z\in X$ and such that $D\cdot l=0$. This would imply
$D\cdot Q=2 D\cdot l=0$, which is clearly impossible. In
conclusion $R_x=X$, so that  the Picard number of $X$ is one by
\cite[IV.3.13.3]{Kollar}. In this case,  clearly
$\Pic(X)=\mathbb{Z}\langle \O(1)\rangle $. Moreover, $Q$ has ample normal bundle, so $2i(X)=-K_X \cdot Q \geq n+1$, giving $i(X)\geq \frac{n+1}2$. 
\vskip4pt
From now on, assume that $[l_1]\not = [l_2]$, so that $V_1\not = V_2$. 
\vskip4pt

\noindent{\it Case} II. Assume that both $V_1$ and $V_2$ are covering. So $a_i\geq 0$, $i=1,2$.

\vskip4pt

Keeping the notation  introduced above,  let us observe that
\begin{equation}\label{stimadim}\begin{split}
\dim(F^x_1)+\dim(F^x_2)&\geq\dim(C^x_1)+\dim(C^x_2)
=a_1+1 + a_2 +1\\&=-K_X\cdot
l_1-1+(-K_X)\cdot l_2-1\\&=-K_X\cdot Q-2=\dim(\mathcal{C}_x)\geq
n-1.\end{split}\end{equation}
Moreover, since $[l_1]\neq[l_2]$, an $R_1$-equivalence class and an $R_2$-equivalence class cannot intersect along a curve, see \cite[IV.3.13.3]{Kollar}. Note that we shall use this systematically in what follows. 

Assume first that $\dim(F^x_1)+\dim(F^x_2) = n$. If $y\in X$ is another general point, using the rational maps $\pi_i$ and the fact that $\dim(F^x_1\cap F^x_2) \leq 0$, we see that $F^x_1\cap F^y_2\neq\emptyset$. In particular, any point of $X$ may be joined to $x$ by a connected chain of rational curves from the families $V_1, V_2$. It follows from \cite[IV.3.13.3]{Kollar} that $N_1(X)_{\mathbb{Q}}= \langle [l_1],[l_2]\rangle$. Moreover, by \cite[Proposition 1]{BCD}, both rays $\mathbb{R}^+[l_i]$ are extremal since $R_i$-equivalence classes are equidimensional, for $i=1,2$. 

 Assume next that $\dim(F^x_1) + \dim(F^x_2) = n - 1$. By (2.1) this gives $a_1 + a_2 + 2 = n - 1$. Denote by $E^x_i$ the (closed) locus of points lying on rational curves from the family $V_j$ meeting $F^x_i$, for $i\neq j$. We have $\dim(E^x_i) = a_1 + a_2 + 2 = n - 1$, $i=1,2$. It follows from the definition of $E^x_i$ that $E^x_i\cdot l_j=0$ for $i\neq j$. We must have $E^x_i\cdot l_i > 0$; otherwise we get $E^x_i\cdot Q=0$, which is absurd since $E^x_i$ is effective and $x$ is general. From this it follows again that any point of $X$ may be joined to $x$ by a connected chain of rational curves from the families $V_1, V_2$, so that $N_1(X)_{\mathbb{Q}}=\langle[l_1],[l_2]\rangle$. We claim that both rays $\mathbb{R}^+[l_1], \mathbb{R}^+[l_2]$ are extremal. Take a curve $C\subset X$ and write its class in $N_1(X)_{\mathbb{Q}}$ as $[C]= a_1[l_1] + a_2[l_2]$, $a_1,a_2\in \mathbb{Q}$. We may assume (up to a permutation of indices) that $a_1\geq 0$. If $E^x_2\cdot C\geq 0$, it follows $a_2\geq 0$ and we are done. If not, we must have $C\subset E^x_2$. But in this case, a variant of ``bend and break" shows that we may write $[C] = \alpha_1[l_1] + \alpha_2[l_2]$, with $\alpha_2\geq 0$. To the best of our knowledge, this last trick was first noticed in \cite[Lemma 1.4.5]{BSW}.

In conclusion, the cone of effective $1$-cycles 
$NE(X)= \overline{NE(X)}=\mathbb{R}^+[l_1]+\mathbb{R}^+[l_2]$ and
 $-K_X$ is ample by the Kleiman criterion. Thus $X$ is a Fano
manifold with
$\Pic(X)\simeq\mathbb{Z}\oplus\mathbb{Z}$. 

There exist two contractions $f_i:X\to Z_i$ of the extremal rays ${\mathbb R} ^+ [l_i]$, where $Z_i$
is a normal projective variety with $\Pic (Z_i)\simeq \mathbb{Z}$ for  $i=1,2$. Write, for simplicity, $F_i=F^x_i$ and $C_i=C_i^x$ for $i=1,2$. 

From \eqref{stimadim} it follows that exactly one of the following subcases (a), (b) or (c) may arise. 
{\small \begin{alignat*}{2}
\mbox{\rm (a) }\,&F_i= C_i\simeq \p^{a_i+1}, i=1,2,\, \dim (F_1)+\dim (F_2)=n, &\quad&\mbox{so } a_1+a_2+2=n,\\
\mbox{\rm (b) }\,&F_i= C_i\simeq \p^{a_i+1},  i=1,2,\,\dim (F_1)+\dim (F_2)=n-1, &&\mbox{so } a_1+a_2+2=n-1,\\
\mbox{\rm (c) }\, &F_1= C_1\simeq \p^{a_1+1}, C_2\subsetneq F_2, \,\dim (F_1)+\dim (F_2)=n,  &&\mbox{so } a_1+a_2+3=n.
 \end{alignat*}}
\noindent{\it Subcase} (a). The linear spaces $F_1,F_2$ intersect transversally in a point. It follows that the map 
$ f_1\times f_2: X \to Z_1\times Z_2 \simeq \p^{a_2+1} \times \p^{a_1+1}$ is finite and birational, hence an isomorphism. This leads to case (iv). 
\vskip4pt

\noindent{\it Subcase} (b).
At least one of the contractions $f_1, f_2$ is equidimensional. Otherwise we may find two fibres $F_1,F_2$ such that $\dim(F_1\cap F_2) \geq 1$. Assume that $f_1$ is equidimensional. From \cite[Lemma\deft2.12]{Fu} or \cite[p.\ 467]{Io1}, it follows that $f_1: X\to Z_1$ is a $\p^{a_1+1}$-bundle. In particular, $Z_1$ is smooth. At this point, we want to apply the characterization of projective space due to Cho--Miyaoka--Shepherd-Barron \cite{CMSB} to see that $Z_1\simeq \p^{a_2+2}$. As the proof of this important criterion is quite involved, we would like to explain how, due to our particular context, the result follows from the simpler argument in \cite{Ke2}. Note that smoothness of $Z_1$ is important for this simplification.  Let $Q$ be a conic from our family $\mathcal {C}_x$. If $Q= l_1\cup l_2$ is reducible, we have $f_1(Q)= f_1(l_2)$ and the restriction $f_1|l_2: l_2\to f_1(l_2)$ is an isomorphism. Let now $Q$ be irreducible and let $\Pi$ denote the plane it spans. We claim that $f_1|Q: Q\to f_1(Q)$ is birational. This is seen by looking at the intersection of the linear spaces $\Pi$ and a fibre $F_1$. If $x\in X$ is a very general point, consider the family of irreducible reduced rational curves $f_{1*}(Q)$, where $Q\in \mathcal{C}_x$. The argument in \cite{Ke2} applies to this family, yielding that $Z_1\simeq \p^{a_2+2}$ and $f_1(Q)$ is a line. 

A general conic $Q$ from our family  has ample normal bundle and $\deg N_{Q|X}= a_1+a_2+2= n-1$, so $Q$ is a quasi-line.
Now fix a general conic $Q_0$ and let $L:= f_1(Q_0)$, a line in $Z_1= \p^{a_2+2}$. Let $T:= f_1^{-1} (L) $. Fix also two general points $t,t'\in T$. Any conic $Q$ passing through $t,t'$ is contained in $T$, because $f_1(Q)$ is a line. The standard    exact sequence of normal bundles
\[0 \to N_{Q|T} \to N_{Q|X}\to N_{T|X\, |Q}\to 0\]
shows that $Q$ is a quasi-line in $T$.   
Therefore, we may apply \cite[Proposition\deft4.1]{IV} to deduce that (up to a Stein factorization) the only non-constant non-finite morphisms defined on $T$ are $f_1|T$ and $\sigma:T\to \p^{a_1+2}$, a blowing-up sending $Q$ to a line of $ \p^{a_1+2}$. Consider now the restriction $r:=f_2|T : T\to Z_2$. As $\dim (T)= a_1+2$, by the above considerations, $r$ is surjective. We claim that $r$ is birational. Otherwise, take two distinct points $t_1, t_2$ in a general fibre, $ F_2\cap T$, of $r$. The line $\langle t_1, t_2\rangle \subset F_2$ meets $T$ in two points, so it is contained in $T$ by the above reasoning. This is absurd since $r$ is generically finite. Moreover, $r$ is not an isomorphism, as $\Pic (Z_2)\simeq \mathbb{Z}$ and $\Pic (T) \simeq \mathbb{Z}\oplus \mathbb{Z} $. Therefore, $r=\sigma$, in particular $Z_2 \simeq \p^{a_1+2}$ and $f_2(Q)$ is a line. 

To conclude subcase (b), we write the hyperplane section of $X$ as $H=H_1+H_2$, with $H_i\cdot l_j= 1 - \delta_{ij}$. We have seen that $f_i^*(D_i) = H_i$, where $D_i$ is the hyperplane divisor of $Z_i$ for $i=1,2$. The map $f:= f_1\times f_2 :X\to \p^{a_2+2} \times \p ^{a_1+2}$ is given by a linear system $\Lambda\subseteq |H|$. Let $X':= f(X)$, which is a divisor on $Z:=   \p^{a_2+2} \times \p ^{a_1+2} $, and note that $f:X\to X'$ is finite and birational by construction. Since $H\cdot l_i = 1$ for $i = 1,2$, the general fiber of $f_i$ is mapped by $f$ onto a hyperplane of $Z_j$, for $i\neq j$. Therefore, the divisor $X'$ has type $(1,1)$ in Pic$(Z)$ and $f$ induces an isomorphism on each fiber of $f_1$. Thus $f:X\to X'$ is an isomorphism and we get case (iii) of the theorem. 
\vskip4pt

\noindent{\it Subcase} (c). In this case, both $f_1$ and $f_2$ are equidimensional and 
$\dim (Z_2)= n-a_2-2=a_1+1=\dim (F_1)$. Using the notation and the reasoning above we get $\dim (T) = a_1+2$, so $f_2|T:T \to Z_2$ cannot be birational. This shows subcase (c) does not happen.  

\vskip4pt

\noindent{\it Case} III. Assume that the family $V_1$ is covering and $V_2$ is not.

 Now we can  also assume that the maximal $x$-covering family $\mathcal{C}_x$ is smooth (otherwise we are in Case II for some different covering family $V_2$), so we may use Theorem \ref {4.5}. 
 In particular, the general conic in $\mathcal{C}_x$ is a quasi-line. 

\vskip4pt

We consider separately the subcases:
\begin{itemize}
\item[III$_1$.] $C^x_1$ has codimension $1$ in $X$.

\item[III$_2$.] $C^x_1$ has codimension $\geq 2$ in $X$.
\end{itemize}

\noindent{\it Case} III$_1$. Reasoning as in Case II, we deduce that either
$X=F^x_1\supsetneq C_1^x$, or $F^x_1=C_1^x$ is a linear subspace in $\p^N$. In the first case, $\Pic(X)=\mathbb{Z}\langle \O(1)\rangle $ and we would fall in Case I, contrary to our assumption.
In the second case, $F^x_1$ has trivial normal bundle, so that  $X\subset\p^N$ is easily seen to 
be a scroll over $\p^1$, say $\pi:X\to \p^1$, having $F^x_1$ as general fiber. We may write $X \simeq \p(\mathcal E)$, with $\mathcal E = \pi_*(\mathcal{O}_X(1))$ a rank-$n$ ample vector bundle on $\p^1$. Since  the general conic $Q \in \mathcal{C}_x$ is a quasi-line,  $-K_X\cdot Q=n+1$. Using the formula for the canonical class of $\p(\mathcal E)$ and the relation $F^x_1\cdot Q=1$, we obtain $\deg(\det\mathcal E)=n+1$. As $\mathcal {E}$ is ample, its splitting-type must be $(1,\ldots,1,2)$ and we get case (ii) with $r=1$ (or $s=n-2$).

\vskip4pt
 
\noindent{\it Case} III$_2$. 

{\it Claim}: If $Q$ is a general member of $\mathcal{C}_x$, $Q$ meets the locus of degenerate conics only at $x$. 
\vskip4pt

Let $D\subset X$ denote the locus of $V_2$. By \cite[II.3.7]{Kollar}, to prove the Claim it is sufficient to show that $D\cap Q=\emptyset$ when  $D$ is a prime divisor. 
In this case   $C_1^x$ is also  irreducible and the locus of degenerate conics is $C^x_1 \cup D$.

The strategy for proving the above Claim is in two steps. First, as in Case II, we show that the rationally connected fibration determined by the covering family $V_1$ is a $\p^ {a_1+1}$-bundle over $\p^{n-a_1-1}$.
Next, we deduce from this that a general conic through $x$ does not intersect $D$, 
proving the Claim. We now proceed with the details.

Let $y$ be a general point of $F^x_1\cap D$. As before, let $C^y_2$ be the cone of lines from $V_2$ passing through $y$. The geometrical description of $D$ yields
\[\dim (C^x_1 \cap D) + \dim (C^y_2) = \dim (D) = n-1.\]
Since \[\dim (C^x_1\cap D) = \dim (C^x_1) -1 =-K_X \cdot l_1-2,\] 
we find    
\[\dim (C^y_2) = n-1 + K_X \cdot l_1+2= -K_X \cdot Q + K_X \cdot l_1= - K_X \cdot l_2 \quad \hbox{for $y\in D$ general.}\]
It follows
\begin{equation}\label{2.2}
\dim (C_2^y) \geq -K_X \cdot l_2\quad \hbox{ for any $y\in D$.}
\end{equation}
By generality of $x$ and $y$ we may assume that 
$(F_1^x \cap D) \cap C^y_2$ meets the smooth locus of $D$.
We deduce 
\begin{equation*}\begin{split}
0&\geq  \dim ((F _1^x \cap D) \cap C_2^y) \geq \dim(F_1^x \cap D) + \dim (C_2^y ) - (n-1) \\& \geq \dim ( C_1^x \cap D) + \dim (C_2^y ) - (n-1) =0.\end{split} 
\end{equation*}
In particular, $\dim (F^x_1 \cap D) = \dim (C^x_1 \cap D)$. Thus, $F^x_1 = C^x_1 = \p^{a_1+1}$, $a_1=\deg (N_{l_1|X})\geq 0$. 
Next we prove that the rationally connected fibration determined by $V_1$ is a $\p^ {a_1+1}$-bundle over some manifold $Z$. Indeed, let $\pi: X^0 \to Z^0$ be our fibration, which is generically a $\p^ {a_1+1}$-bundle by the above. We show that $X^0=X$ and $\pi$ is equidimensional. Using the notation in \cite[Proposition 1]{BCD}, assume that $B\not =\emptyset$. 
We claim that $B\cap D \not= \emptyset$. Indeed, if $b\in B\setminus D$, by definition of $B$ there exists a line $l_1\in V_1$ contained in $B$ and containing $b$. Since $D\cdot l_1>0$,  the last claim follows. By \cite[Proposition 1]{BCD}, there is some equivalence class $\overline{F_1}$ contained in $B$, passing through a point $b\in B\cap D$ and having $\dim ( \overline{F_1} ) > a_1+1$. We deduce, using \eqref{2.2} 
\[\dim (\overline {F_1} \cap C_2^b)\geq \dim (\overline{F_1}) + \dim (C_2^b)- n > a_1 +1 - K_X \cdot l_2 -n =0.\] 
This would imply $[l_1]=[l_2]$, a contradiction. So, $B=\emptyset$ and $\pi $ is equidimensional. 
As in Case II, subcase (b), it follows that $\pi: X\to Z$ is a $\p^ {a_1+1}$-bundle and $Z\simeq\p^{n-a_1-1}$. Moreover, a general conic
through $x$ is sent isomorphically by $\pi$ onto a line in $Z$.

Let us fix $Q$, a general conic through $x$. Set $L:=\pi(Q)\subset\p^{n-a_1-1}$ and let $T:=\pi^{-1}(L)\to L$ be the associated $\p^{a_1+1}$-bundle over $\p^1$.
The conic $Q$  is contained in $T$. Reasoning as in Case II, subcase (b), it is immediate to see that
$T$ is a CC-manifold, belonging to the class described in Case III$_1$ above. In particular there is a map $\sigma: T\to \p^{a_1+2}$ which is a blowing-up along a linear $\p^{a_1}$.
General conics in $T$ are identified via $\sigma$ to lines in $\p^{a_1+2}$ and $D\cap T$ is the exceptional divisor of $\sigma$. 
Thus $Q\cap D=Q\cap (D\cap T)=\emptyset$ and the Claim is finally proved.
\vskip4pt

 By  Theorem \ref{4.5}(ii), there
exists a birational map $\phi:\p^n\map X\subset\p^N$ given by a linear
system  $\Lambda \subseteq |\O_{\p^n}(2)|$, inducing an isomorphism between a neighbourhood of a line in $\p^n$ and a neighbourhood of a conic $Q\in \mathcal{C}_x$. 
 In particular $\phi$ is an isomorphism outside the base locus of $\Lambda$, which we denote by $Y$. 
  We claim that, if $Y\neq\emptyset$,
$Y$ is a linear space $\p^s$ with  $0\leq s\leq n-3$.

Suppose there exist two distinct points $y_1, y_2\in Y$. 
The linear system of quadric hypersurfaces $\Lambda$ cannot separate points along the line $\langle y_1,y_2\rangle $.  Therefore $\langle y_1, y_2\rangle \subseteq Y$. The infinitesimal version of the same argument shows $Y$ is reduced, so it is a linear space.
 Thus we are in one of the cases described in (ii).

\vskip4pt

Finally, if all conics in $\mathcal{C}_x$ are irreducible, we get case (i) of the theorem. This follows directly from Theorem~\ref{4.5}(ii).  
\end{proof}

Conversely, from the recent paper \cite{BH}, we see that the condition $i\geq \frac{n+1}2$ in Theorem \deft\ref{delta1} is optimal.

\begin{Proposition}\label{2.4}Let $X\subset \p ^{n+r}$ be a smooth non-degenerate complete intersection of dimension $n$ and multi-degree $(d_1, d_2, \ldots, d_r)$. If we assume that $n\geq 2\sum_1^r d_i-2r-1$, then $X$ is a CC-manifold.  
\end{Proposition}
\begin{proof} Note that the inequality in the statement means that $X$ is Fano, of index at least $\frac{n+1}2$. Replacing $X$ by a general hyperplane section, we may assume that equality holds. Therefore, we may apply \cite[Corollary 1.2]{BH} to conclude.
\end{proof}

The next result, essentially due to Hwang--Kebekus \cite[Theorem\deft3.14]{HK}, shows that any Fano manifold of high index is conic-connected. Note that we slightly improve the bound on the index given in \cite{HK}. The first part of the proposition is well known.

\begin{Proposition}[cf.\ \cite{HK}]\label{largeindex}
 Let  $X\subset\p^N$ be a Fano manifold with $\Pic(X)\simeq \Z \langle H\rangle$ and $-K_X = i(X)H$, $H$ being the hyperplane section and $i(X)$ the index~of~$X$.
\begin{enumerate}
\item[(i)] If $i(X)>\frac{n+1}{2}$, then $X\subset\p^N$ is ruled by
lines and for general $x\in X$ the Hilbert scheme of lines through
$x$, $Y_x\subset \p( T ^*_x(X))=\p^{n-1}$, is   smooth.
If $i(X)\geq\frac{n+3}{2}$, $Y_x$ is also irreducible.
Let $SY_x$ be the secant variety of $Y_x$ in $\p^{n-1}$.
\item[(ii)] If $i(X)\geq\frac{n+3}{2}$ and  $SY_x=\p^{n-1}$, then
$X\subset\p^N$ is a CC-manifold.
\item[(iii)] If
$i(X)>\frac{2n}{3}$, then $X\subset\p^N$ is a CC-manifold. 
\end{enumerate}
\end{Proposition}
\begin{proof} By a theorem of Mori \cite{Mori} the variety $X\subset\p^N$ is
ruled by a family of rational curves $\mathcal{L}$ such that for
$L\in\mathcal{L}$ we have $-K_X\cdot L\leq n+1$. 
It follows that $X\subset\p^N$ is
ruled by lines, that is $H\cdot L=1$. As is well known, see \cite[II.3.11.5]{Kollar} for a more general result, the Hilbert scheme of lines passing
through $x$ is  smooth equidimensional and may be
identified with a subscheme
$Y_x\subset\p( T 
^*_x(X))=\p^{n-1}$ of dimension
$i(X)-2$. Thus, if $i(X)-2\geq \frac{n-1}2$, $Y_x\subset\p^{n-1}$ is irreducible.

Part (ii) follows from \cite[Theorem 3.14]{HK}. 

To prove part (iii), first observe that $i(X) \geq \frac{n+3}2$, unless $n\leq 6$. If $n\leq 6$, $i(X)>\frac{2n}3$ gives $i(X)\geq n-1$, so the conclusion follows by the classification of del Pezzo manifolds (see \cite{Fujita-book}). To conclude, using part (ii), it is enough to see that $SY_x=\p^{n-1}$.

 By \cite[Theorem
2.5]{Hwang}, the variety $Y_x\subset\p^{n-1}$ is non-degenerate and
by hypothesis $n-1<\frac{3i(X)-6}{2}+2=\frac{3\dim(Y_x)}{2}+2$, so
that $SY_x=\p^{n-1}$ by Zak Linear Normality Theorem; see
\cite[V.1.13]{Zak}. 
\end{proof}

\section*{Acknowledgements}

Both authors are grateful to the organizers of the Conference ``Projective varieties with unexpected properties", that took place in Siena, between 8--13 of June 2004. Our collaboration started while taking part in this very pleasant and fruitful mathematical event. 

We thank Tommaso de Fernex and Ciro Ciliberto for a critical reading of previous  versions of our paper.

\end{document}